\documentclass[11pt]{article}
\usepackage[a4paper,margin=25mm]{geometry}
\usepackage{amsmath,amssymb,amsthm,mathtools}
\usepackage[hidelinks]{hyperref}
\usepackage{comment}
\usepackage{tikz}
\usepackage{placeins}
\usepackage{caption}
\usepackage{xcolor}
\usetikzlibrary{patterns}

\numberwithin{equation}{section}
\allowdisplaybreaks[1]

\theoremstyle{plain}
\newtheorem{theorem}{Theorem}[section]
\newtheorem{lemma}[theorem]{Lemma}
\newtheorem{proposition}[theorem]{Proposition}
\newtheorem{corollary}[theorem]{Corollary}
\theoremstyle{remark}
\newtheorem{remark}[theorem]{Remark}

\newcommand{\R}{\mathbb{R}}
\newcommand{\N}{\mathbb{N}}
\newcommand{\dd}{\,\mathrm{d}}
\newcommand{\tm}{\tau_-}

\newcommand{\tP}{\tau_{\rm P}}
\newcommand{\norm}[1]{\left\|#1\right\|}
\newcommand{\la}{\langle}
\newcommand{\ra}{\rangle}
\newcommand{\acp}{\alpha_{\rm C}}

\title{Critical exponent for Nakao's problem in three space dimensions}
\author{{\large Vladimir Georgiev}\\[3mm]
	    \normalsize Department of Mathematics, 
        \\\normalsize University of Pisa, Largo Bruno Pontecorvo 5, 56127, Pisa, Italy
        \\\normalsize (e-mail: georgiev@dm.unipi.it)
        \\\normalsize (ORCID: \url{https://orcid.org/0000-0001-6796-7644})\\[4mm]
{\large Kosuke Kita\footnote{Corresponding author}}\\[3mm]
        \normalsize Department of Mathematics, Hokkaido University,
        \\\normalsize Kita-10, Nishi-8, Kita-ku, Sapporo, Hokkaido, 060-0810, Japan
        \\\normalsize (e-mail: kosuke.kita@math.sci.hokudai.ac.jp)
        \\\normalsize (ORCID: \url{https://orcid.org/0000-0001-7099-3253})}
\date{}

\begin{document}
\maketitle

\begin{abstract}
We study a system consisting of a semilinear damped wave equation coupled to a semilinear wave equation in three space dimensions.
We determine the sharp condition on the exponents of nonlinear terms that separates global existence for sufficiently small initial data from finite-time blow-up even for arbitrarily small data satisfying suitable positivity assumptions.
To prove global existence, we combine a weighted pointwise estimate for the damped wave equation, which we established in earlier work, with a new one for the wave equation.
The new pointwise estimate for the wave equation uses a combination of a parabolic weight and a light-cone weight.
For the proof of blow-up, pointwise lower bounds are compared with integral upper bounds obtained by test functions.
In the main estimates, explicit physical-space representations reveal how forcing near the wave cone influences the damped component on parabolic spatial scales and, in turn, the wave component inside the cone.
\end{abstract}



\section{Introduction}

\subsection{Nakao's problem and related topics}

This paper concerns the following Cauchy problem
\begin{equation}\label{eq:nakao}
\left\{
\begin{aligned}
 &\partial_t^2 u-\Delta u+\partial_t u=|v|^p, && t>0,\ x\in\R^3,\\
 &\partial_t^2 v-\Delta v=|u|^q, && t>0,\ x\in\R^3,\\
 &(u,\partial_t u,v,\partial_t v)(0,x)=(u_0,u_1,v_0,v_1)(x), && x\in\R^3,
\end{aligned}
\right.
\end{equation}
where \(p,q>1\).  
This system couples a semilinear damped wave equation with a semilinear wave equation.  
The problem of finding the threshold between global existence and finite-time blow-up for \eqref{eq:nakao} was proposed by Professor Mitsuhiro Nakao as a way of comparing the nonlinear effect of damping with that of wave propagation; the system is therefore often called Nakao's problem (see \cite{Nishihara-Wakasugi2015}).

The damped wave equation is hyperbolic and has finite propagation speed, but its large-time behavior contains a heat-like part.  
More precisely, for suitably localized data, solutions of the linear damped wave equation approach solutions of the heat equation with a suitable initial datum.  
This is the diffusion phenomenon; see, for example, \cite{Matsumura,Chill-Haraux,Ikehata-Nishihara,Nishihara1997,Yang-Milani}.  
The linear propagators of the two components in \eqref{eq:nakao} therefore have different large-time behavior: the damped component has heat-like decay, whereas the undamped component retains the light-cone structure of the wave equation.

We seek a quantity \(\acp=\acp(p,q)\) whose sign characterizes the threshold between global existence and finite-time blow-up for \eqref{eq:nakao} in the following sense:
\begin{itemize}
    \item[(i)] If \(\acp<0\), then there exists a unique global solution to \eqref{eq:nakao} for sufficiently small initial data satisfying suitable regularity and spatial decay assumptions.
    \item[(ii)] If \(\acp\ge0\), then finite-time blow-up occurs even for arbitrarily small initial data satisfying suitable positivity assumptions.
\end{itemize}
The precise data and solution classes will be specified in the main results below.
Once such a quantity is identified, the set \(\acp(p,q)=0\) is called the \emph{critical curve}, and the pairs \((p,q)\) on this curve are called \emph{pairs of critical exponents} (for short, \emph{critical exponents} or \emph{critical pairs}).

For the scalar heat equation with nonnegative data
\begin{equation*}
    \left\{
    \begin{aligned}
        &\partial_t u-\Delta u=u^p, &&\qquad t>0,~x\in \R^n,\\
        &u(0)=u_0,
    \end{aligned}
    \right.
\end{equation*}
the critical exponent is \(p_{\rm F}(n)=1+2/n\).
Fujita \cite{Fujita} proved finite-time blow-up for nontrivial data when \(1<p<p_{\rm F}(n)\) and global existence for sufficiently small, suitably decaying data when \(p>p_{\rm F}(n)\).
The critical case also belongs to the blow-up range, as shown by Hayakawa \cite{Hayakawa} in dimensions one and two and by Kobayashi, Sirao and Tanaka \cite{Kobayashi-Sirao-Tanaka} in arbitrary dimensions; see also Weissler \cite{Weissler}.

For the scalar wave equation in dimensions \(n\ge2\)
\begin{equation*}
    \left\{
    \begin{aligned}
        &\partial_t^2 u-\Delta u=|u|^p, &&\qquad t>0,~x\in \R^n,\\
        &u(0)=u_0,~\partial_t u(0)=u_1,
    \end{aligned}
    \right.
\end{equation*}
the Strauss exponent \(p_{\rm S}(n)\) is the positive root of \((n-1)p^2-(n+1)p-2=0\).
In three dimensions, John \cite{John} proved global existence for small, regular, compactly supported data above \(p_{\rm S}(3)=1+\sqrt2\) and blow-up below it under suitable positivity assumptions.
Strauss \cite{Strauss} proposed the threshold in general dimensions. For \(n\ge3\), the first author, Lindblad and Sogge \cite{Georgiev-Lindblad-Sogge} proved global existence for small, smooth, compactly supported data when \(p_{\rm S}(n)<p\le(n+3)/(n-1)\); see \cite{Sogge} for an account of the scalar theory.
For suitable compactly supported data, critical blow-up was proved by Schaeffer \cite{Schaeffer} for \(n=2,3\) and by Yordanov and Zhang \cite{Yordanov-Zhang} for \(n\ge4\), where Takamura and Wakasa \cite{Takamura-Wakasa} obtained sharp upper lifespan bounds.
Lifespan estimates in low dimensions were obtained by Lindblad \cite{Lindblad} and Zhou \cite{Zhou1,Zhou2,Zhou3}. In the higher-dimensional subcritical case, Di Pomponio and the first author \cite{DiPomponio-Georgiev} derived lower bounds, and Takamura \cite{Takamura} obtained upper bounds.

For the scalar damped wave equation
\begin{equation*}
    \left\{
    \begin{aligned}
        &\partial_t^2 u-\Delta u + \partial_t u=|u|^p, &&\qquad t>0,~x\in \R^n,\\
        &u(0)=u_0,~\partial_t u(0)=u_1,
    \end{aligned}
    \right.
\end{equation*}
the critical exponent for localized data is again \(p_{\rm F}(n)\).
For compactly supported energy data, Todorova and Yordanov \cite{Todorova-Yordanov} proved small-data global existence above \(p_{\rm F}(n)\) (with \(p\le n/(n-2)\) when \(n\ge3\)) and blow-up below it under an integral positivity assumption. Zhang \cite{Zhang} established nonexistence at the endpoint under a positivity condition by a test function argument.
Ikehata and Tanizawa \cite{Ikehata-Tanizawa} replaced compact support in the global result by spatial decay. In two dimensions, Hosono and Ogawa \cite{Hosono-Ogawa} proved global existence and diffusion-type asymptotics for the related nonlinearity \(|u|^{p-1}u\), \(p>2\), for sufficiently small data in suitable Besov spaces intersected with \(L^1\).
In three dimensions, Nishihara \cite{Nishihara2003} used linear \(L^p\)-\(L^q\) estimates to prove global existence and optimal decay of small weak solutions for powers above \(p_{\rm F}(3)\).
Upper lifespan bounds below and at \(p_{\rm F}(n)\) were obtained by Li and Zhou \cite{Li-Zhou} for \(n=1,2\) and by Nishihara \cite[Theorem~1.4]{Nishihara2003-lifespan} for \(n=3\). Ikeda and Wakasugi \cite{Ikeda-Wakasugi} extended the subcritical upper bounds to arbitrary dimensions. For critical lifespan estimates, see Ikeda and Ogawa \cite{Ikeda-Ogawa} and, for the sharp exponential upper bound in dimensions \(n\ge4\), Lai and Zhou \cite{Lai-Zhou}.

For two weakly coupled heat equations 
\begin{equation*}
    \left\{
    \begin{aligned}
        &\partial_t u-\Delta u =v^p, &&\qquad t>0,~x\in \R^n,\\
        &\partial_t v-\Delta v =u^q, &&\qquad t>0,~x\in \R^n,
    \end{aligned}
    \right.
\end{equation*}
with nonnegative initial data, the critical curve is given by
\(\alpha_{\rm DW}=0\), where
\begin{equation*}
    \alpha_{\rm DW}
    :=\max\left\{\frac{p+1}{pq-1},\frac{q+1}{pq-1}\right\}-\frac n2.
\end{equation*}
Escobedo and Herrero \cite{Escobedo-Herrero} proved nonexistence of nontrivial global solutions for \(\alpha_{\rm DW}\ge0\) and global existence for sufficiently small, suitably decaying data for \(\alpha_{\rm DW}<0\).
The same index governs the system of two damped wave equations with powers \(|v|^p,|u|^q\). Sun and Wang \cite{Sun-Wang} proved small-data global existence for \(\alpha_{\rm DW}<0\) in dimensions one and three, and nonexistence for \(\alpha_{\rm DW}\ge0\) in arbitrary dimensions under integral positivity assumptions.
Narazaki \cite{Narazaki} included dimension two and treated slowly decaying data. Nishihara and Wakasugi \cite{Nishihara-Wakasugi} proved small-data global existence in arbitrary dimensions under the power and localization restrictions of their weighted energy method.

For two weakly coupled wave equations
\begin{equation*}
    \left\{
    \begin{aligned}
        &\partial_t^2 u-\Delta u =|v|^p, &&\qquad t>0,~x\in \R^n,\\
        &\partial_t^2 v-\Delta v =|u|^q, &&\qquad t>0,~x\in \R^n,
    \end{aligned}
    \right.
\end{equation*}
the corresponding index for \(n\ge2\) is
\begin{equation*}
        \alpha_{\rm W}:=
        \max\left\{\frac{p+2+q^{-1}}{pq-1},\frac{q+2+p^{-1}}{pq-1}\right\}-\frac{n-1}{2}.
\end{equation*}
Note that, on the diagonal \(p=q\), the curves \(\alpha_{\rm DW}=0\) and \(\alpha_{\rm W}=0\) recover the Fujita and Strauss exponents, respectively.
Del Santo, the first author and Mitidieri \cite{DelSanto-Georgiev-Mitidieri} proved small-data global existence and blow-up results on the two sides of \(\alpha_{\rm W}=0\) under additional assumptions on the initial data and the exponents \(p\) and \(q\).
Kubo and Ohta \cite{Kubo-Ohta} proved critical blow-up and upper lifespan bounds for suitable data in dimensions two and three. For lifespan estimates including the critical case in three dimensions, see Agemi, Kurokawa and Takamura \cite{Agemi-Kurokawa-Takamura}.
In dimensions \(n\ge4\), see \cite{Georgiev-Takamura-Zhou} for subcritical lifespan estimates and \cite{Kurokawa-Takamura-Wakasa} for critical blow-up and upper lifespan bounds.
Ikeda, Sobajima and Wakasa \cite{Ikeda-Sobajima-Wakasa} obtained lifespan estimates by a test function method that uses integral positivity rather than pointwise nonnegativity of the data.

For Nakao's problem, Wakasugi \cite{Wakasugi}, Chen and Reissig \cite{Chen-Reissig}, and the second author and Kusaba \cite{Kita-Kusaba} proved nonexistence of global solutions under their respective assumptions on the data, powers and solution class. 
In three dimensions, their strict conditions give
\begin{equation}\label{eq:known-blowup-pq}
\alpha_{\rm BU}(p,q;3)
:=\max\left\{
\frac{p+1}{pq-1}-\frac32,
\frac{q+1}{pq-1}-\frac32,
\frac{p^{-1}+2}{pq-1}-1
\right\}>0.
\end{equation}
The energy-space results summarized in \cite{Kita-Kusaba} impose \(1<p,q\le3\) in three dimensions. Some endpoint cases are also covered in the cited works.

In \cite[Corollaries~1.3--1.5]{Georgiev-Kita-weighted}, we used weighted pointwise estimates for the Duhamel term to construct global mild solutions for small data in three dimensions when
\begin{equation*}
 q>\frac{4p+2}{p^2}\quad (1<p\le3),
 \qquad
 q>1+\frac5{3p}\quad (p>3).
\end{equation*}
For sufficiently small, spatially integrable Sobolev data, Chen \cite{Chen-global} proved global existence in three dimensions for \(q>2\) and \(p>\max\{2,1+10/(3q)\}\), as well as scattering of the wave component in \(H^2\times H^1\).
This exponent range lies within the one above, but the data classes and conclusions differ.
Moreover, he treated the two-dimensional case.

The previous existence and blow-up results leave a gap. We show that, in three dimensions, the global-existence range for sufficiently small weighted data and the blow-up range for the positive data specified below meet at the curve \(\alpha_{\rm C}(p,q;3)=0\), where
\begin{equation}\label{eq:alpha-CP}
        \alpha_{\rm C}(p,q;3)
        :=\frac{p+1}{pq-1}-\min\left\{\frac p2,\frac32\right\}.
\end{equation}

The two exponents in the minimum have different origins. 
The exponent \(3/2\) is the decay rate of the three-dimensional heat kernel, whereas \(p/2\) reflects how the forcing \(|v|^p\) near the light cone contributes to the damped component in the region \(r\lesssim\sqrt{t}\).
To capture both effects, the principal pointwise estimates and blow-up arguments use explicit physical-space propagators in the \((t,r)\)-plane, \(r=|x|\). 
The new pointwise estimate for the Duhamel term of the wave equation gives decay in \(|t-r|\) from a parabolic decay estimate on its source \(|u|^q\).
Combined with the Duhamel estimate of the damped wave equation \cite{Georgiev-Kita-weighted}, it yields global existence throughout the range \(\alpha_{\rm C}<0\).

The role of this spatial dependence is clearest in the blow-up proof on the critical segment \(1<p<3\). 
In the critical case \(q=(3p+2)/p^2\), \(1<p<3\), the earlier blow-up index satisfies \(\alpha_{\rm BU}<0\), yet we prove finite-time blow-up even for arbitrarily small data satisfying
the positivity assumptions below.
A lower bound for the free wave near \(r=t\) first yields a lower bound for the damped component on \(r\le\sqrt{t}/4\). 
Inserting this bound back into the wave equation gives a lower bound for \(|v|^p\) with the factor \((t-r)^{-1}\) inside the light cone. 
Integration in \(r\) then produces the logarithm needed to contradict the upper bound from the parabolic test function.

\subsection{Main results}

We write \(L^\infty_{t,x}:=L^\infty((0,\infty)\times\R^3)\).
For \(M\in\N\) and \(\xi\in\R^M\), we set \(\la\xi\ra:=1+|\xi|\).
For \(t\ge0\) and \(x\in\R^3\), define
\begin{equation*}
        \tm=\tm(t,x):=\la t-|x|\ra,
        \qquad
        \tP=\tP(t,x):=\la t+|x|^2\ra.
\end{equation*}
For nonnegative quantities \(A\) and \(B\), we write \(A\lesssim B\) if \(A\le CB\) for a constant \(C>0\) independent of the variables under consideration, and write \(A\simeq B\) if both \(A\lesssim B\) and \(B\lesssim A\) hold.  

Let \(W(t)\) be the wave propagator, that is, \(W(t)f\) is the solution of
\begin{equation*}
    \left\{
    \begin{aligned}
        &\partial_t^2 w-\Delta w=0, &&t>0,~x\in\R^3\\
        &w(0)=0,\quad\partial_t w(0)=f.
    \end{aligned}
    \right.
\end{equation*}
Kirchhoff's formula gives
\begin{equation}\label{eq:W-explicit}
        W(t)f(x)=\frac{t}{4\pi}\int_{\mathbb S^2}f(x+t\omega)\,\dd\omega.
\end{equation}
Let \(D(t)\) be the corresponding propagator for
\begin{equation*}
    \left\{
    \begin{aligned}
        &\partial_t^2 w-\Delta w + \partial_t w=0, &&t>0,~x\in\R^3\\
        &w(0)=0,\quad\partial_t w(0)=f.
    \end{aligned}
    \right.
\end{equation*}
In three space dimensions, this propagator has the following representation (see, e.g., \cite{Nishihara2003}):
\begin{equation}\label{eq:D-decomposition}
        D(t)=e^{-t/2}W(t)+S(t),
\end{equation}
where
\begin{equation}\label{eq:S-explicit}
        S(t)f(x)
        :=\frac{e^{-t/2}}{4\pi}
        \int_0^t\int_{\mathbb S^2}
        \frac{I_1\left(\frac12\sqrt{t^2-\rho^2}\right)}
        {2\sqrt{t^2-\rho^2}}
        \rho^2 f(x+\rho\omega)\,\dd\omega\,\dd\rho.
\end{equation}
Here \(I_1\) is the modified Bessel function of the first kind of order one,
\begin{equation*}
        I_1(s)=\sum_{k=0}^{\infty}
        \frac{1}{k!(k+1)!}\left(\frac{s}{2}\right)^{2k+1}.
\end{equation*}
By Duhamel's principle, \eqref{eq:nakao} can be written as the following integral equation.
\begin{equation}\label{eq:mild-system}
\left\{
\begin{aligned}
        u&=K_D[u_0,u_1]+L_D[|v|^p],\\
        v&=K_W[v_0,v_1]+L_W[|u|^q],
\end{aligned}
\right.
\end{equation}
where, for \(\varphi_0=\varphi_0(x)\), \(\varphi_1=\varphi_1(x)\) and \(F=F(t,x)\),
\begin{equation*}
\begin{aligned}
        K_D[\varphi_0,\varphi_1](t)&:=\partial_tD(t)\varphi_0+D(t)(\varphi_0+\varphi_1),\\
        K_W[\varphi_0,\varphi_1](t)&:=\partial_tW(t)\varphi_0+W(t)\varphi_1,
\end{aligned}
\end{equation*}
and
\begin{equation*}
        L_D[F](t):=\int_0^tD(t-s)F(s)\,\dd s,
        \qquad
        L_W[F](t):=\int_0^tW(t-s)F(s)\,\dd s.
\end{equation*}

We first state a global existence result for \eqref{eq:nakao}.
Since our proof relies on space-time weighted \(L^\infty\)-estimates, we introduce suitable function spaces. 
For \(\alpha, \beta>0\), define
\begin{align*}
        &X_\alpha:=\{F\in L^\infty_{t,x} ~\!;~\! \norm{F}_{X_\alpha}<\infty\},\qquad
        \norm{F}_{X_\alpha}:=
        \left\|\tP^\alpha F(t,x)\right\|_{L^\infty_{t,x}},\\[2mm]
        &Y_\beta:=\{F\in L^\infty_{t,x} ~\!;~\! \norm{F}_{Y_\beta}<\infty\},\qquad
        \norm{F}_{Y_\beta}
        :=\left\|\tP^{1/2}\tm^\beta F(t,x)\right\|_{L^\infty_{t,x}}.
\end{align*}
Equipped with these norms, \(X_\alpha\) and \(Y_\beta\) are Banach spaces.
For \((\varphi_0,\varphi_1)\in W^{1,\infty}(\R^3)\times L^\infty(\R^3)\), we also set
\begin{align*}
\mathcal N_D^\alpha(\varphi_0,\varphi_1)
&:=\|\la x\ra^{2\alpha}\varphi_0\|_{L^\infty}
  +\|\la x\ra^{2\alpha}\nabla \varphi_0\|_{L^\infty}
  +\|\la x\ra^{2\alpha}\varphi_1\|_{L^\infty},\\[2mm]
\mathcal N_W^\beta(\varphi_0,\varphi_1)
&:=\|\la x\ra^{1+\beta}\varphi_0\|_{L^\infty}
  +\|\la x\ra^{2+\beta}\nabla \varphi_0\|_{L^\infty}
  +\|\la x\ra^{2+\beta}\varphi_1\|_{L^\infty}.
\end{align*}
A pair \((u,v)\in X_\alpha\times Y_\beta\) is called a global mild solution of \eqref{eq:nakao} if \eqref{eq:mild-system} holds almost everywhere in \((0,\infty)\times\R^3\).

\begin{theorem}\label{thm:global}
Let \(p,q>1\), and assume
\begin{equation}\label{eq:global-index-condition}
        \alpha_{\rm C}(p,q;3)<0.
\end{equation}
Choose \(\alpha,\beta>0\) so that
\begin{equation}\label{eq:global-weight-choice}
\alpha<\min\left\{\frac p2,\frac32\right\},
\qquad
\frac{\alpha+1}{p}-\frac12<\beta<\alpha q-\frac32.
\end{equation}
Assume \((u_0,u_1), (v_0,v_1)\in W^{1,\infty}(\R^3)\times L^\infty(\R^3)\) and set \(\mathcal N_0:=\mathcal N_D^\alpha(u_0,u_1)+\mathcal N_W^\beta(v_0,v_1)\).
There exists a constant \(\varepsilon_0>0\) such that, if \(\mathcal N_0\le\varepsilon_0\), then \eqref{eq:nakao} has a unique global mild solution \((u,v)\) satisfying
\begin{equation}\label{eq:global-bound}
        \norm{u}_{X_\alpha}+\norm{v}_{Y_\beta}
        \lesssim \mathcal N_0.
\end{equation}
\end{theorem}

\begin{remark}
If a global mild solution in \(X_\alpha\times Y_\beta\) exists for given initial data, then it is unique in this class, regardless of the smallness of the data.
Since \(D(t)\) and \(W(t)\) preserve nonnegativity and satisfy \(D(t)1=1-e^{-t}\) and \(W(t)1=t\), their \(L^\infty\) operator norms are bounded on every finite time interval.
Then, uniqueness follows by applying Gronwall's inequality to the difference of the integral equations on each such interval.
\end{remark}

\begin{remark}
Condition \eqref{eq:global-index-condition} is equivalent to
\begin{equation*}
        q>\frac{3p+2}{p^2}
        \quad (1<p\le3),
        \qquad
        q>\frac{2p+5}{3p}
        \quad (p>3).
\end{equation*}
For \(p\ge5\), the second inequality follows from \(q>1\).  Moreover, when \(1<p\le3\), the first inequality implies \(pq>3\), and hence \(3/(2q)<p/2\).  When \(p>3\), the inequality \(q>1\) gives \(3/(2q)<3/2\).  Thus one can choose
\[
\max\left\{\frac{3}{2q},\frac{p+1}{pq-1}\right\}
<\alpha<\min\left\{\frac p2,\frac32\right\}.
\]
The interval for \(\beta\) in \eqref{eq:global-weight-choice} then contains positive numbers.
\end{remark}

To state our blow-up result precisely, we first need to clarify the notion of blow-up.
We therefore begin by establishing the local existence of solutions via a standard argument.
The following result ensures the existence of a maximal solution, which will be used to formulate our blow-up theorem.
Its proof is given in Appendix~\ref{app:local-sobolev}.
\begin{proposition}\label{prop:local-sobolev}
Let \(p,q>1\).  For every \((u_0,u_1),(v_0,v_1)\in H^2(\R^3)\times H^1(\R^3)\), \eqref{eq:nakao} has a unique maximal solution \((u,v)\) on \([0,T_m)\), where \(0<T_m\le\infty\), such that
\begin{equation}\label{eq:local-sobolev-regularity-main}
        (u,v)\in
        \left(\bigcap_{j=0}^2 C^j\bigl([0,T_m);H^{2-j}(\R^3)\bigr)\right)^2.
\end{equation}
Moreover, either \(T_m=\infty\), or 
\begin{equation}\label{eq:local-sobolev-Linfty-alt-main}
        T_m<\infty,\qquad
        \limsup_{t\uparrow T_m}
        \left(\norm{u(t)}_{L^\infty}+\norm{v(t)}_{L^\infty}\right)=\infty.
\end{equation}
\end{proposition}

The blow-up alternative obtained by standard arguments for evolution equations takes the form
\begin{equation*}
    \lim_{t\uparrow T_m}
        \left(\norm{u(t)}_{H^2}+\norm{v(t)}_{H^2}+\norm{\partial_t u(t)}_{H^1}+\norm{\partial_t v(t)}_{H^1}\right)=\infty
\end{equation*}
whenever \(T_m<\infty\).
In the present setting, the weaker \(L^\infty\) norm suffices for the blow-up alternative, although \(H^2(\R^3)\) embeds into \(L^\infty(\R^3)\). Consequently, Theorem~\ref{thm:global} gives the following result.

\begin{corollary}\label{cor:compatibility}
Let \((u_0,u_1),(v_0,v_1)\in H^2(\R^3)\times H^1(\R^3)\).
In addition, assume that the hypotheses of Theorem~\ref{thm:global} hold.
There exists a constant \(\varepsilon_0>0\) such that, if \(\mathcal N_0\le\varepsilon_0\), then \(T_m=\infty\) and \eqref{eq:global-bound} hold.
\end{corollary}

For the blow-up result, we consider profiles \(\phi,\psi\) such that
\begin{equation}\label{eq:phi-psi-class}
        \phi,\psi\in C_c^\infty(\R^3),
        \qquad
        \phi,\psi\ge0
\end{equation}
satisfying
\begin{equation}\label{eq:phi-psi-pos}
        \phi(x)>0,
        \qquad
        \psi(x)>0
        \qquad (|x|\le1),
\end{equation}
\begin{equation}\label{eq:phi-psi-support}
        \operatorname{supp}\phi\cup\operatorname{supp}\psi
        \subset\{x\in\R^3:|x|\le2\}.
\end{equation}

\begin{theorem}\label{thm:blowup}
Let \(p,q>1\), and assume
\begin{equation}\label{eq:blowup-index-condition}
        \alpha_{\rm C}(p,q;3)\ge0.
\end{equation}
Let \((\phi,\psi)\) satisfy \eqref{eq:phi-psi-class}--\eqref{eq:phi-psi-support} and assume 
\begin{equation}\label{eq:data}
        u_0=v_0=0,
        \qquad
        u_1=\phi,
        \qquad
        v_1=\psi.
\end{equation}
Let \((u,v)\) be the maximal solution as in Proposition~\ref{prop:local-sobolev}.
Then \(T_m<\infty\); that is, the solution blows up in finite time in the sense of \eqref{eq:local-sobolev-Linfty-alt-main}.
\end{theorem}

\begin{remark}\label{rem:alphaCP-blowup-equivalence}
For \(p,q>1\), condition \eqref{eq:blowup-index-condition} is equivalent to one of the following alternatives:
\begin{align*}
        &1<p<3,
        ~ 1<q\le\frac{3p+2}{p^2},\\
        &3\le p<5,
        ~1<q\le\frac{2p+5}{3p}.
\end{align*}
For \(p\ge5\), there is no \(q>1\) satisfying \eqref{eq:blowup-index-condition}.  
\end{remark}

To compare the existence of global solutions (Theorem \ref{thm:global} and Corollary \ref{cor:compatibility}) and the nonexistence of global solutions (Theorem \ref{thm:blowup}) on the same class of initial data, take \(\phi,\psi\) satisfying \eqref{eq:phi-psi-class}--\eqref{eq:phi-psi-support} and set \eqref{eq:data}.
Then
\begin{equation*}
        \mathcal N_0
        =\mathcal N_D^\alpha(0,\phi)
        +\mathcal N_W^\beta(0,\psi).
\end{equation*}
These profiles can be chosen so that \(\mathcal N_0\) is arbitrarily small.
Thus sufficiently small \(\mathcal N_0\) gives a global solution when
\(\alpha_{\rm C}(p,q;3)<0\), whereas finite-time blow-up occurs without
any size restriction when \(\alpha_{\rm C}(p,q;3)\ge0\).

\begin{figure}[htbp]
\centering
\begin{tikzpicture}[x=1.35cm,y=1.50cm]
\pgfmathsetmacro{\pzero}{(3+sqrt(41))/8}
\pgfmathsetmacro{\pone}{1+sqrt(10)/2}
\pgfmathsetmacro{\pA}{5/3}
\pgfmathsetmacro{\pB}{1+sqrt(2)}
\pgfmathsetmacro{\pmax}{8.1}
\pgfmathsetmacro{\qmax}{5.3}

\draw[->] (0,0) -- (8.4,0) node[right] {$p$};
\draw[->] (0,0) -- (0,5.45) node[above] {$q$};
\node[below left] at (0,0) {$O$};

\draw[densely dotted] (1,0) -- (1,5.25);
\draw[densely dotted] (3,0) -- (3,5.25);
\draw[densely dotted] (5,0) -- (5,5.25);
\draw[densely dotted] (0,1) -- (8.2,1);

\draw (1,0.04) -- (1,-0.04) node[below] {$1$};
\draw (3,0.04) -- (3,-0.04) node[below] {$3$};
\draw (5,0.04) -- (5,-0.04) node[below] {$5$};
\draw (0.04,1) -- (-0.04,1) node[left] {$1$};
\draw (0.04,5) -- (-0.04,5) node[left] {$5$};

\begin{scope}
\clip (1,1) rectangle (3,\qmax);
\path[pattern=north east lines,pattern color=black]
  (1.005,1)
  -- plot[domain=1.005:3,samples=80,smooth,variable=\p]
       ({\p},{(3*\p+2)/(\p*\p)})
  -- (3,1)
  -- cycle;
\end{scope}
\begin{scope}
\clip (3,1) rectangle (5,\qmax);
\path[pattern=north east lines,pattern color=black]
  (3,1)
  -- plot[domain=3:5,samples=60,smooth,variable=\p]
       ({\p},{(2*\p+5)/(3*\p)})
  -- (5,1)
  -- cycle;
\end{scope}

\begin{scope}
\clip (1,1) rectangle (3,\qmax);
\path[pattern=vertical lines,pattern color=black]
  (1.005,\qmax)
  -- (3,\qmax)
  -- plot[domain=3:1.005,samples=80,smooth,variable=\p]
       ({\p},{(3*\p+2)/(\p*\p)})
  -- cycle;
\end{scope}
\begin{scope}
\clip (3,1) rectangle (5,\qmax);
\path[pattern=vertical lines,pattern color=black]
  (3,\qmax)
  -- (5,\qmax)
  -- plot[domain=5:3,samples=60,smooth,variable=\p]
       ({\p},{(2*\p+5)/(3*\p)})
  -- cycle;
\end{scope}
\path[pattern=vertical lines,pattern color=black]
  (5,1) rectangle (\pmax,\qmax);

\draw[thin]
  plot[domain=1.005:\pzero,samples=60,smooth,variable=\p]
       ({\p},{5/(3*\p-2)});
\draw[thin]
  plot[domain=\pzero:\pone,samples=60,smooth,variable=\p]
       ({\p},{(3*\p+1)/(\p*\p)});
\draw[thin]
  plot[domain=\pone:5,samples=60,smooth,variable=\p]
       ({\p},{(2*\p+5)/(3*\p)});

\draw[densely dashed,very thick]
  plot[domain=1.005:\pA,samples=60,smooth,variable=\p]
       ({\p},{5/(3*\p-2)});
\draw[densely dashed,very thick]
  plot[domain=\pA:5,samples=80,smooth,variable=\p]
       ({\p},{(2*\p+5)/(3*\p)});

\draw[dash dot,very thick]
  plot[domain=1.68:\pB,samples=60,smooth,variable=\p]
       ({\p},{(3*\p+1)/(\p*(\p-1))});
\draw[dash dot,very thick]
  plot[domain=\pB:\pmax,samples=80,smooth,variable=\p]
       ({\p},{(\p+3+sqrt(\p*\p+10*\p+9))/(2*\p)});

\draw[very thick]
  plot[domain=1.005:3,samples=80,smooth,variable=\p]
       ({\p},{(3*\p+2)/(\p*\p)});
\draw[very thick]
  plot[domain=3:5,samples=60,smooth,variable=\p]
       ({\p},{(2*\p+5)/(3*\p)});

\node[fill=white,inner sep=1.5pt] at (1.95,1.31)
{\bfseries\small Blow-up};
\node[fill=white,inner sep=1.5pt] at (4.15,3.65)
{\bfseries\small Global existence};

\begin{scope}[shift={(4.50,5.00)}]
\draw[very thick] (0,0) -- (0.72,0);
\node[fill=white,right] at (0.82,0) {\small $\alpha_{\rm C}=0$};

\draw[thin] (0,-0.26) -- (0.72,-0.26);
\node[fill=white,right] at (0.82,-0.26) {\small $\alpha_{\rm BU}=0$};

\draw[densely dashed,very thick] (0,-0.54) -- (0.72,-0.54);
\node[fill=white,right] at (0.82,-0.54) {\small $\alpha_{\rm DW}=0$};

\draw[dash dot,very thick] (0,-0.82) -- (0.72,-0.82);
\node[fill=white,right] at (0.82,-0.82) {\small $\alpha_{\rm W}=0$};
\end{scope}
\end{tikzpicture}
\caption{The global-existence and finite-time blow-up regions in Theorems~\ref{thm:global} and~\ref{thm:blowup}}
\label{fig:pq-plane-global}
\end{figure}
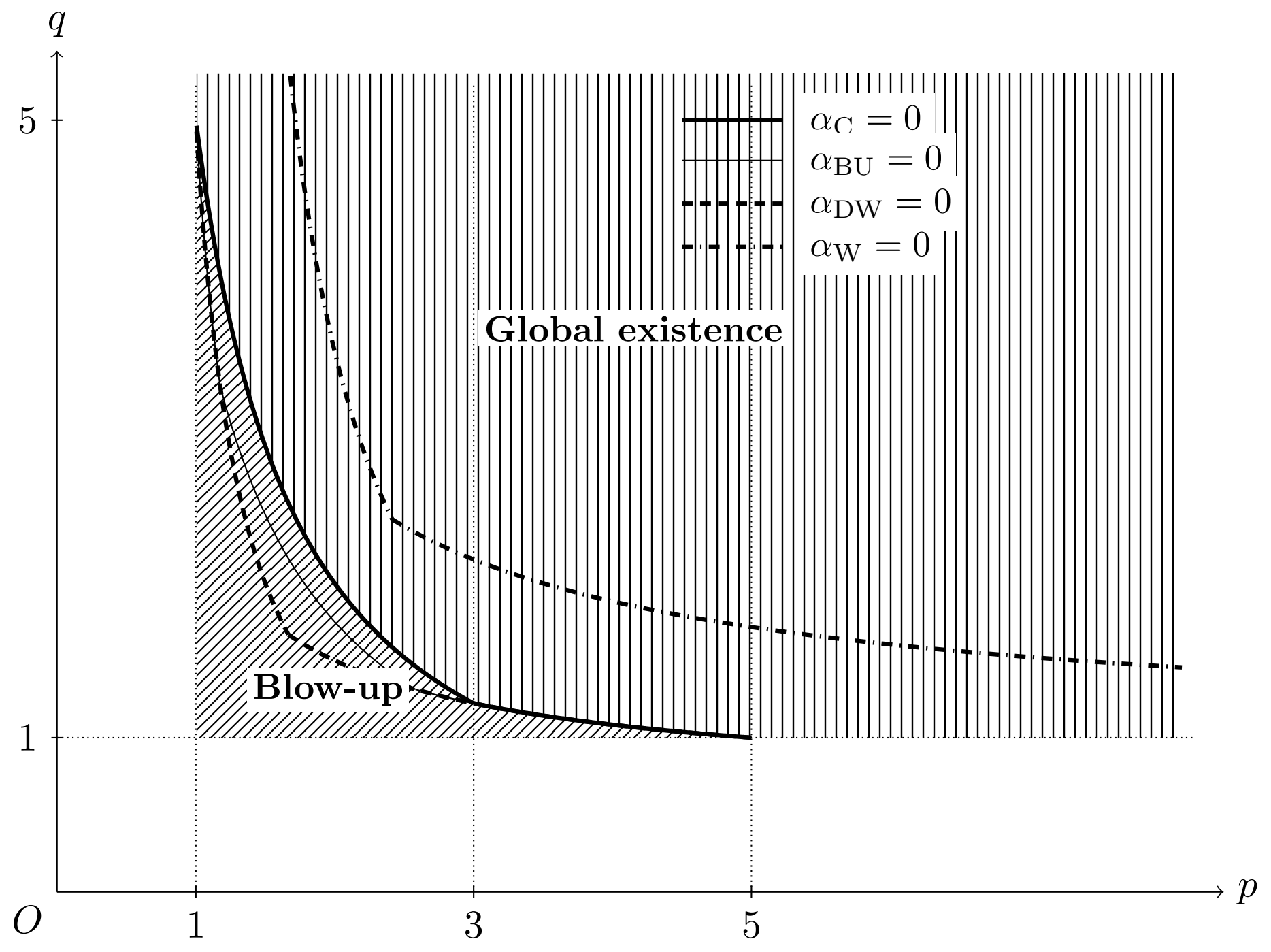
\FloatBarrier

Figure~\ref{fig:pq-plane-global} shows the global-existence and finite-time blow-up regions established in Theorems~\ref{thm:global} and~\ref{thm:blowup}.  
The boundary of the strict nonexistence region for global weak solutions in \eqref{eq:known-blowup-pq} and the critical curves for the damped wave--damped wave and wave--wave systems are included for comparison.

The rest of the paper is organized as follows.  Section~2 presents the weighted estimates and proves Theorem~\ref{thm:global} by a contraction argument in \(X_\alpha\times Y_\beta\).  Section~3 proves Theorem~\ref{thm:blowup}, first in the strict case and then on the two critical segments.  The appendix proves Proposition~\ref{prop:local-sobolev} and Corollary~\ref{cor:compatibility}.

\section{Global existence}

\subsection{Weighted estimates for the linear damped wave equation}

We first recall the weighted estimates for the homogeneous linear
solutions.  The notation for the propagators and the data norms is the
one introduced in the statement of the main result.

\begin{lemma}\label{lem:weighted-homogeneous}
Let \(0<\alpha<3/2\).  Then
\begin{equation*}
        \norm{K_D[u_0,u_1]}_{X_\alpha}
        \lesssim \mathcal N_D^\alpha(u_0,u_1).
\end{equation*}
\end{lemma}

This estimate is stated in \cite[Proposition~2.3]{Georgiev-Kita-weighted};
see also \cite{Narazaki-Nishihara}.
The following Duhamel estimate for the damped wave equation is
\cite[Theorem~1.1]{Georgiev-Kita-weighted}.

\begin{proposition}\label{prop:damped-duhamel}
Let \(0<\alpha\le3/2\), \(\tilde{\alpha}>\alpha\), and \(\gamma>1\).  Then
\begin{equation}\label{eq:damped-duhamel}
\norm{L_D[F]}_{X_\alpha}
\lesssim
\left\|\tP^{\tilde{\alpha}}\tm^\gamma F\right\|_{L^\infty_{t,x}}
\end{equation}
whenever the right-hand side is finite.
\end{proposition}

\subsection{New weighted estimates for the linear wave equation}

We shall use the standard spherical reduction formula: if \(x\in\R^3\),
\(r=|x|>0\), \(\rho>0\), and \(h:[0,\infty)\to[0,\infty]\) is
measurable, then
\begin{equation}\label{eq:spherical-reduction-global}
\int_{\mathbb S^2}h(|x+\rho\omega|)\,\dd\omega
=\frac{2\pi}{r\rho}
\int_{|r-\rho|}^{r+\rho}\lambda h(\lambda)\,\dd\lambda.
\end{equation}
We next state the elementary integral estimate needed below.

\begin{lemma}\label{lem:radial-integrals-global}
Let \(\beta>0\).  
For \(a,\sigma\ge0\),
\begin{equation}\label{eq:one-dimensional-global}
\int_0^a
\la\lambda+(\lambda-\sigma)^2\ra^{-\beta-1/2}\,\dd\lambda
\le C_\beta \la\sigma\ra^{-\beta},
\qquad\quad
C_\beta=2^{\beta+1/2}\frac{\sqrt{\pi} \Gamma\left(\beta\right)}{\Gamma\left(\beta+\frac{1}{2}\right)}.
\end{equation}
\end{lemma}

\begin{proof}
We can see that
\begin{equation*}
1+\lambda+(\lambda-\sigma)^2
\ge 
\frac12 (1+\sigma+(\lambda-\sigma)^2).
\end{equation*}
Hence, 
\begin{equation*}
\int_0^a
\la\lambda+(\lambda-\sigma)^2\ra^{-\beta-1/2}\,\dd\lambda
\le 2^{\beta+1/2} \int_\R(1+\sigma+(\lambda-\sigma)^2)^{-\beta-1/2}\,\dd \lambda.
\end{equation*}
Setting \(\lambda-\sigma=\sqrt{1+\sigma}~\!\tau\), we have
\begin{equation*}
    \int_\R(1+\sigma+(\lambda-\sigma)^2)^{-\beta-1/2}\,\dd \lambda
    =
    \la \sigma\ra^{-\beta} \int_\R (1+\tau^2)^{-\beta-1/2} \dd\tau. 
\end{equation*}
The change of variables \(\eta=\tau^2\) gives
\begin{align*}
    \int_\R (1+\tau^2)^{-\beta-1/2} \dd\tau
    &=
    2\int_0^\infty (1+\tau^2)^{-\beta-1/2} \dd\tau\\[2mm]
    &=\int_0^\infty \eta^{-1/2} (1+\eta)^{-\beta-1/2}\dd\eta\\[2mm]
    &=B\left(\frac{1}{2},\beta\right)=\frac{\sqrt{\pi} \Gamma\left(\beta\right)}{\Gamma\left(\beta+\frac{1}{2}\right)}.
\end{align*}
This completes the proof.
\end{proof}

\begin{lemma}\label{lem:wave-homogeneous-parabolic}
Let \(\beta>0\).  Then
\begin{equation*}
        \norm{K_W[v_0,v_1]}_{Y_\beta}
        \lesssim \mathcal N_W^\beta(v_0,v_1).
\end{equation*}
\end{lemma}
This lemma follows from \cite[Lemma~3.1, (3.1)]{Katayama-Yokoyama}
together with the inequality \(\tP^{1/2}\le\la t+r\ra\).

\begin{proposition}\label{prop:wave-duhamel-parabolic}
Let \(\beta>0\).  Then
\begin{equation*}
\norm{L_W[F]}_{Y_\beta}
\lesssim
\left\|\tP^{\beta+3/2}F\right\|_{L^\infty_{t,x}}
\end{equation*}
whenever the right-hand side is finite.
\end{proposition}

\begin{proof}
Kirchhoff's formula gives
\begin{align*}
|L_W[F](t,x)|
\le\frac{1}{4\pi}\int_0^t(t-s)
\int_{\mathbb S^2}
\la s+|x+(t-s)\omega|^2\ra^{-\beta-3/2}
\,\dd\omega\,\dd s ~\!\left\|\tP^{\beta+3/2}F\right\|_{L^\infty_{t,x}}.
\end{align*}
Set \(r=|x|\) and
\begin{equation*}
    I(t,x):= \int_0^t(t-s)
\int_{\mathbb S^2}
\la s+|x+(t-s)\omega|^2\ra^{-\beta-3/2}
\,\dd\omega\,\dd s
\end{equation*}
It is sufficient to show 
\begin{equation}
    \label{eq:wave_esti}
    I(t,x) \lesssim \tP^{-1/2} \tm^{-\beta}.
\end{equation}
We first suppose that \(0\le t\le16\).  
Since \(t-s\le16\), we find \(\la x+(t-s)\omega\ra\simeq\la r\ra\) for \(0\le s\le t\).
Indeed, this equivalence follows from
\begin{equation*}
    \la x+(t-s)\omega\ra \le 1+r+t-s\le 17 \la r\ra
\end{equation*}
and
\begin{equation*}
    \la r\ra \le 1+ |x+(t-s)\omega| +t-s\le 17 \la x+(t-s)\omega\ra.
\end{equation*}
Hence, 
\begin{align*}
    I(t,x)\lesssim \la r\ra^{-2\beta-3}.
\end{align*}
Since \(\tP^{1/2}\tm^{\beta}\lesssim\la r\ra^{1+\beta}\), we obtain
\begin{equation*}
    I(t,x)\lesssim \tP^{-3/2} \tm^{-2\beta}\le \tP^{-1/2} \tm^{-\beta}.
\end{equation*}

It remains to consider \(t>16\).
We further divide the argument into two cases.

\noindent\underline{(i) Case of \(0\le r<\sqrt{t}\)}:
Setting \(\rho=t-s\), we have
\begin{equation*}
    I(t,x)=\int_0^t \rho \int_{\mathbb{S}^2}\la t-\rho +|x+\rho\omega|^2\ra ^{-\beta-3/2}\dd\omega\dd\rho.
\end{equation*}
If \(0\le \rho \le2\sqrt t\), then \(t-\rho\ge t-2\sqrt t\ge t/2\) due to \(t>16\).
On the other hand, it holds that \(|x+\rho\omega|\ge \rho- r\ge \rho-\sqrt t\ge \rho/2\)  provided that \(2\sqrt t\le\rho \le t\).
Therefore, 
\begin{align*}
    I(t,x)
    &\lesssim 
    \la t\ra^{-\beta-3/2} \int_0^{2\sqrt t} \rho \dd\rho
    +
    \int_{2\sqrt t}^{t} \rho \la \rho^2\ra^{-\beta-3/2}\dd\rho\\[2mm]
    &\lesssim
    \la t\ra^{-\beta-1/2}
    \lesssim
    \tP^{-1/2}\tm^{-\beta},
\end{align*}
where we used \(\tP\lesssim\la t\ra\) and \(\tm\simeq\la t\ra\) in this case.

\noindent\underline{(ii) Case of \(\sqrt{t}\le r\)}:
By virtue of \eqref{eq:spherical-reduction-global}, it holds that
\begin{equation*}
    I(t,x)\lesssim \frac{1}{r}\int_0^t\int_{|r-t+s|}^{r+t-s}\lambda\la s+\lambda^2\ra^{-\beta-3/2}\dd\lambda\dd s.
\end{equation*}
Here, we can verify
\begin{align*}
    \int_a^\infty \lambda \la s+\lambda^2\ra^{-\beta-3/2}\dd\lambda
    &=\frac{1}{2}\int_{s+a^2}^\infty \la\eta\ra^{-\beta-3/2}\dd\eta\\[2mm]
    &=\frac{1}{2\beta+1}\la s+a^2\ra^{-\beta-1/2}
    \qquad
    (s,a\ge 0),
\end{align*}
which implies
\begin{align*}
    I(t,x)
    &\lesssim 
    \frac{1}{r}\int_0^t \la s+ (r-t+s)^2\ra^{-\beta-1/2}\dd s\\[2mm]
    &\le 
    \frac{1}{r}\int_0^t \la s+ (s-|t-r|)^2\ra^{-\beta-1/2}\dd s,
\end{align*}
where we used \((r-t+s)^2\ge (s-|t-r|)^2\) for \(0\le s\le t\).
Applying Lemma \ref{lem:radial-integrals-global} to the right hand side with \(a=t\) and \(\sigma=|t-r|\) of the above inequality, we obtain
\begin{equation*}
    I(t,x)\lesssim r^{-1}\tm^{-\beta}.
\end{equation*}
Note that \(\tP^{1/2}\simeq r\) follows from \(t>16\) and \(\sqrt t\le r\).
Thus, the above estimate is equivalent to \eqref{eq:wave_esti} in this case.
\end{proof}

\subsection{Existence of global mild solutions}

\begin{proof}[Proof of Theorem~\ref{thm:global}]
By \eqref{eq:global-weight-choice}, choose \(\tilde{\alpha}>\alpha\) and \(\gamma>1\)
sufficiently close to \(\alpha\) and \(1\), respectively, so that
\begin{equation}\label{eq:global-parameter-relations}
\tilde{\alpha}<\frac p2,
\qquad
\tilde{\alpha}+\gamma<p\left(\beta+\frac12\right),
\qquad
\beta+\frac32<\alpha q.
\end{equation}

On \(X_\alpha\times Y_\beta\), use the product norm
\begin{equation*}
\norm{(u,v)}_{X_\alpha\times Y_\beta}
:=\norm{u}_{X_\alpha}+\norm{v}_{Y_\beta}.
\end{equation*}
Lemmas~\ref{lem:weighted-homogeneous} and~\ref{lem:wave-homogeneous-parabolic} give
\begin{equation*}
\norm{K_D[u_0,u_1]}_{X_\alpha}
+\norm{K_W[v_0,v_1]}_{Y_\beta}
\lesssim \mathcal N_0.
\end{equation*}

For \(v\in Y_\beta\),
\begin{equation*}
|v(t,x)|
\le\norm{v}_{Y_\beta}\tP^{-1/2}\tm^{-\beta}.
\end{equation*}
By \eqref{eq:global-parameter-relations}, \(p/2-\tilde{\alpha}>0\) and \(\gamma-p\beta<p/2-\tilde{\alpha}\).
Since \(1\le\tm\lesssim\tP\),
\begin{equation*}
\tP^{\tilde{\alpha}-p/2}\tm^{\gamma-p\beta}
\le \left(\frac{\tm}{\tP}\right)^{p/2-\tilde{\alpha}}
\lesssim 1.
\end{equation*}
Hence
\begin{equation*}
\left\|\tP^{\tilde{\alpha}}\tm^\gamma|v|^p\right\|_{L^\infty_{t,x}}
\lesssim \norm{v}_{Y_\beta}^p.
\end{equation*}
Therefore, Proposition~\ref{prop:damped-duhamel} yields
\begin{equation*}
\norm{L_D[|v|^p]}_{X_\alpha}
\lesssim \norm{v}_{Y_\beta}^p.
\end{equation*}

Similarly, for \(u\in X_\alpha\),
\begin{equation*}
\tP^{\beta+3/2}|u|^q
\le\norm{u}_{X_\alpha}^q
\tP^{\beta+3/2-\alpha q}.
\end{equation*}
The exponent on the right is negative by
\eqref{eq:global-parameter-relations}.  Proposition~\ref{prop:wave-duhamel-parabolic} gives
\begin{equation*}
\norm{L_W[|u|^q]}_{Y_\beta}
\lesssim\norm{u}_{X_\alpha}^q.
\end{equation*}

Define \(\Phi\) on \(X_\alpha\times Y_\beta\) by
\begin{equation*}
\Phi(u,v):=
\left(K_D[u_0,u_1]+L_D[|v|^p],
K_W[v_0,v_1]+L_W[|u|^q]\right).
\end{equation*}
The estimates above imply
\begin{equation}\label{eq:global-map-bound}
\norm{\Phi(u,v)}_{X_\alpha\times Y_\beta}
\le C_1\mathcal N_0
+C_1\norm{v}_{Y_\beta}^p
+C_1\norm{u}_{X_\alpha}^q
\end{equation}
for some \(C_1>0\).
The mean value inequality, followed by the same weighted estimates,
also implies that there exists \(C_2>0\) such that
\begin{align}
&\norm{\Phi(u,v)-\Phi(\widetilde u,\widetilde v)}_{X_\alpha\times Y_\beta}
\notag\\
&\quad\le C_2\left(\norm{v}_{Y_\beta}^{p-1}
+\norm{\widetilde v}_{Y_\beta}^{p-1}\right)
\norm{v-\widetilde v}_{Y_\beta}
+C_2\left(\norm{u}_{X_\alpha}^{q-1}
+\norm{\widetilde u}_{X_\alpha}^{q-1}\right)
\norm{u-\widetilde u}_{X_\alpha}.
\label{eq:global-map-difference}
\end{align}

Let \(C_0:=\max\{1, C_1,C_2\}\) and set
\begin{equation*}
R:=4C_0\mathcal N_0.
\end{equation*}
Since \(p,q>1\), there is \(\varepsilon_0>0\) such that, whenever
\(\mathcal N_0\le\varepsilon_0\),
\begin{equation*}
C_0(R^p+R^q)\le2C_0\mathcal N_0,
\qquad
2C_0(R^{p-1}+R^{q-1})\le\frac12.
\end{equation*}
It follows that \(\Phi\) maps the closed ball
\begin{equation*}
\left\{(u,v)\in X_\alpha\times Y_\beta:
\norm{(u,v)}_{X_\alpha\times Y_\beta}\le R\right\}
\end{equation*}
into itself and is a contraction there.  Banach's fixed point theorem
therefore gives a global mild solution \((u,v)\) of
\eqref{eq:mild-system} in this ball.  In particular,
\begin{equation*}
\norm{u}_{X_\alpha}+\norm{v}_{Y_\beta}
\le4C_0\mathcal N_0,
\end{equation*}
which proves \eqref{eq:global-bound}.
This completes the proof.
\end{proof}

\section{Finite-time blow-up}

Throughout this section, \((u,v)\) denotes the maximal solution associated with the profiles \eqref{eq:phi-psi-class}--\eqref{eq:phi-psi-support} and the data \eqref{eq:data}, and \([0,T_m)\) is its maximal existence interval.  
By \eqref{eq:local-sobolev-regularity-main} and the embedding \(H^2(\R^3)\hookrightarrow C_b(\R^3)\), we use the continuous representatives of \(u\) and \(v\).  
The Duhamel identities then agree with the pointwise integral formulas.  
We write
\begin{equation*}
        B_R:=\{x\in\R^3:|x|\le R\}.
\end{equation*}
Since the propagators \(D(t)\) and \(W(t)\) preserve nonnegativity, \(D(t)\phi,W(t)\psi\ge0\).  
The mild equations read
\begin{equation*}
        u=D(t)\phi+L_D[|v|^p],
        \qquad
        v=W(t)\psi+L_W[|u|^q].
\end{equation*}
Consequently,
\begin{equation*}
        u\ge D(t)\phi\ge0,
        \qquad
        v\ge W(t)\psi\ge0,
\end{equation*}
and hence
\begin{align}
        v(t,x)&\ge L_W[u^q](t,x),
        \label{eq:trigger-wave-positivity}\\
        u(t,x)&\ge L_D[v^p](t,x).
        \label{eq:trigger-damped-positivity}
\end{align}
All time intervals below are understood to be intersected with \([0,T_m)\).  
Unless otherwise stated, the constants in this section may also depend on the fixed profiles \(\phi\) and \(\psi\).

\subsection{Integral upper bounds from test functions}

Choose a smooth function \(\eta\) on \([0,\infty)\), with \(0\le\eta\le1\), such that
\[
\eta=1\ \hbox{on }[0,1/2],\qquad \eta=0\ \hbox{on }[1,\infty).
\]
Since \(\eta\) is constant near zero, the function \(x\mapsto\eta(|x|)\) is smooth on \(\R^3\).
Let \(p'=p/(p-1)\), \(q'=q/(q-1)\), and fix an integer \(k\ge2\max\{p',q'\}\).  
For \(1\le T<T_m\), set
\begin{equation*}
        \Psi_T(t,x):=
        \left[\eta\left(\frac tT\right)
        \eta\left(\frac{|x|}{\sqrt T}\right)\right]^k,
\end{equation*}
and
\begin{equation*}
        I_T:=\int_0^T\int_{\R^3}v(t,x)^p\Psi_T(t,x)\,\dd x\,\dd t,
        \qquad
        J_T:=\int_0^T\int_{\R^3}u(t,x)^q\Psi_T(t,x)\,\dd x\,\dd t.
\end{equation*}

\begin{lemma}\label{lem:test-function-upper}
With implicit constants independent of \(T\) and the profiles \(\phi,\psi\), we have
\begin{align}
        I_T+\int_{\R^3}\phi(x)\Psi_T(0,x)\,\dd x
        &\lesssim T^{-1+5/(2q')}J_T^{1/q},
        \label{eq:test-damped}\\
        J_T+\int_{\R^3}\psi(x)\Psi_T(0,x)\,\dd x
        &\lesssim T^{-1+5/(2p')}I_T^{1/p}.
        \label{eq:test-wave}
\end{align}
In particular,
\begin{equation}\label{eq:test-upper}
        I_T\lesssim T^{3/2-(p+1)/(pq-1)},
        \qquad
        J_T\lesssim T^{3/2-(q+1)/(pq-1)}.
\end{equation}
\end{lemma}

\begin{proof}
By the construction of the maximal Sobolev solution and the regularity \eqref{eq:local-sobolev-regularity-main}, both equations hold in \(L^2(\R^3)\) for \(0\le t<T_m\).
Pairing the first equation with \(\Psi_T\) in \(L^2\), we may use an integration by parts in time and space to obtain
\begin{equation*}
\begin{aligned}
        I_T
        &=\int_0^T\int_{\R^3}
        (\partial_t^2u-\Delta u+\partial_tu)\Psi_T\,\dd x\,\dd t\\
        &=\left[\int_{\R^3}
        \bigl(\partial_tu\,\Psi_T-u\,\partial_t\Psi_T+u\Psi_T\bigr)
        \,\dd x\right]_0^T
        +\int_0^T\int_{\R^3}
        u(\partial_t^2-\Delta-\partial_t)\Psi_T\,\dd x\,\dd t.
\end{aligned}
\end{equation*}
The terminal boundary terms vanish since
\(\Psi_T(T,\cdot)=\partial_t\Psi_T(T,\cdot)=0\).
Using \(u(0)=v(0)=0\),
\(\partial_tu(0)=\phi\), and \(\partial_tv(0)=\psi\), and
applying the same calculation to the second equation without the
damping term, we obtain
\begin{align*}
        I_T+\int_{\R^3}\phi(x)\Psi_T(0,x)\,\dd x
        &=\int_0^T\int_{\R^3}
        u(\partial_t^2-\Delta-\partial_t)\Psi_T\,\dd x\,\dd t,\\
        J_T+\int_{\R^3}\psi(x)\Psi_T(0,x)\,\dd x
        &=\int_0^T\int_{\R^3}
        v(\partial_t^2-\Delta)\Psi_T\,\dd x\,\dd t.
\end{align*}

We argue as in \cite[Section~2]{Wakasugi}.
Since $T\ge1$ and $\eta$ is constant near zero, a direct calculation gives
\[
|\partial_t\Psi_T|+|\partial_t^2\Psi_T|+|\Delta\Psi_T|
\lesssim T^{-1}\Psi_T^{1-2/k}.
\]
Since \(\operatorname{supp}\Psi_T\subset [0,T]\times B_{\sqrt T}\), we have \(|\operatorname{supp}\Psi_T|\lesssim T^{5/2}\).
By virtue of $k\ge2q'$ and $0\le\Psi_T\le1$, we have $\Psi_T^{1-2/k}\le\Psi_T^{1/q}$. 
Therefore, H\"older's inequality gives
\begin{align*}
I_T+\int_{\R^3}\phi(x)\Psi_T(0,x)\,\dd x
&\lesssim T^{-1}\int_{\operatorname{supp}\Psi_T}
       |u|\Psi_T^{1/q}\,\dd x\,\dd t\\
&\lesssim T^{-1+5/(2q')}J_T^{1/q}.
\end{align*}
Thus we have \eqref{eq:test-damped}.
Likewise, $k\ge2p'$ implies \eqref{eq:test-wave}.
Dropping the nonnegative initial terms and substituting \eqref{eq:test-wave} into \eqref{eq:test-damped}, we obtain
\[
I_T^{1-1/(pq)}\lesssim T^{-1+5/(2q')+(-1+5/(2p'))/q}.
\]
The exponent of $T$ after division by $1-1/(pq)$ equals
\[
\frac{pq}{pq-1}\left(-1+\frac{5}{2q'}
+\frac1q\left(-1+\frac{5}{2p'}\right)\right)
=\frac32-\frac{p+1}{pq-1}.
\]
Similarly, substituting \eqref{eq:test-damped} into \eqref{eq:test-wave} gives the upper bound for \(J_T\) in \eqref{eq:test-upper}.
Note that this proof also shows that the implicit constants are independent of the profiles $\phi,\psi$.
\end{proof}

\subsection{Pointwise lower bounds}

For \(t>0\) and \(|x|<t\), define the positive interior kernel of the propagator of the damped wave equation by
\begin{equation*}
        \mathcal K_D(t,x):=
        \frac{e^{-t/2}}{8\pi}
        \frac{I_1\!\left(\frac12\sqrt{t^2-|x|^2}\right)}
        {\sqrt{t^2-|x|^2}}.
\end{equation*}
Thus, by \eqref{eq:D-decomposition}--\eqref{eq:S-explicit},
\begin{equation*}
        D(t)f(x)\ge S(t)f(x)
        =\int_{|x-y|<t}\mathcal K_D(t,x-y)f(y)\,\dd y
\end{equation*}
whenever \(f\ge0\).

\begin{lemma}\label{lem:damped-kernel}
The following estimate holds.
\begin{equation}\label{eq:damped-lower}
        \mathcal K_D(t,x)\gtrsim t^{-3/2}
\end{equation}
whenever \(t\ge3\) and \(|x|\le\sqrt t\).
\end{lemma}

\begin{proof}
If \(t\ge3\) and \(|x|\le\sqrt t\), then
\begin{equation*}
        0\le t-\sqrt{t^2-|x|^2}
        =\frac{|x|^2}{t+\sqrt{t^2-|x|^2}}
        \le1.
\end{equation*}
Moreover,
\begin{equation*}
        1<\frac{\sqrt6}{2}\le\frac{t}{\sqrt6}
        \le\frac12\sqrt{t^2-|x|^2}\le\frac t2.
\end{equation*}
The standard lower bound
\begin{equation*}
        I_1(R)\gtrsim R^{-1/2}e^R
        \qquad (R\ge1)
\end{equation*}
therefore gives
\begin{align*}
        \mathcal K_D(t,x)
        &\gtrsim t^{-3/2}
        \exp\!\left(-\frac{t-\sqrt{t^2-|x|^2}}2\right)\\
        &\gtrsim t^{-3/2}.
\end{align*}
This proves \eqref{eq:damped-lower}.
\end{proof}

The positive datum \(\psi\) first creates a lower bound near the light cone.

\begin{lemma}\label{lem:wave-cone-lower}
The following estimate holds.
\begin{equation*}
        v(t,x)\gtrsim  t^{-1}
\end{equation*}
for \(1\le t<T_m\) and \(\bigl||x|-t\bigr|\le1/4\).
\end{lemma}

\begin{proof}
Let \(r=|x|\) and \(m=\min_{|y|\le1}\psi(y)>0\).
Since \(\psi\ge m\mathbf{1}_{B_1}\) and \([1/4,1]\subset[|r-t|,r+t]\) for \(t\ge1\) and \(|r-t|\le1/4\), \eqref{eq:W-explicit} and \eqref{eq:spherical-reduction-global} give
\[
    v(t,x)\ge W(t)\psi(x)
    \ge \frac{m}{2r}\int_{1/4}^{1}\lambda\,\dd\lambda
    = \frac{15m}{64r}
    \gtrsim t^{-1},
\]
where \(r\le t+1/4\le5t/4\).
\end{proof}

This lower bound near the light cone (Lemma \ref{lem:wave-cone-lower}) yields a lower bound on parabolic regions.

\begin{lemma}\label{lem:cone-seed}
The following estimate holds.
\begin{equation}\label{eq:seed}
        u(t,x)\gtrsim t^{-\min\{p/2,3/2\}}
\end{equation}
for \(64\le t<T_m\) and \(|x|\le\sqrt t/4\).
\end{lemma}

\begin{proof}
Set
\begin{equation*}
        \mathcal A_s:=\{y\in\R^3:\bigl||y|-s\bigr|\le1/4\}.
\end{equation*}
Lemma~\ref{lem:wave-cone-lower} gives
\begin{equation*}
        v(s,y)^p\gtrsim s^{-p}
        \qquad (s\ge1,\ y\in\mathcal A_s),
\end{equation*}
and
\begin{equation*}
        |\mathcal A_s|
        =\frac{4\pi}{3}\left[\left(s+\frac14\right)^3
        -\left(s-\frac14\right)^3\right]
        =2\pi s^2+\frac{\pi}{24}\ge2\pi s^2.
\end{equation*}
For \(t\ge64\), \(|x|\le\sqrt t/4\),
\(1\le s\le\sqrt t/4\), and \(y\in\mathcal A_s\), we have
\begin{align*}
        t-s&\ge\frac t2\ge32,\\
        |x-y|&\le\frac{\sqrt t}{2}+\frac14
        \le\sqrt{\frac t2}\le\sqrt{t-s}.
\end{align*}
Therefore, applying Lemma~\ref{lem:damped-kernel} to the Duhamel term, we obtain
\begin{align*}
        u(t,x)
        &\gtrsim t^{-3/2}
        \int_1^{\sqrt t/4}s^{-p}|\mathcal A_s|\,\dd s\\
        &\gtrsim t^{-3/2}
        \int_1^{\sqrt t/4}s^{2-p}\,\dd s.
\end{align*}
Here, recall that all source times lie below \(t<T_m\).  
Since \(t\ge64\), we have \(\sqrt t/8\ge1\) and \(\sqrt t/4\ge2\).  
Thus, for \(1<p<3\),
\begin{equation*}
        \int_1^{\sqrt t/4}s^{2-p}\,\dd s
        \ge\int_{\sqrt t/8}^{\sqrt t/4}s^{2-p}\,\dd s
        =\frac{4^{p-3}-8^{p-3}}{3-p}\,t^{(3-p)/2}.
\end{equation*}
If \(p\ge3\), the same integral is at least \(\int_1^2s^{2-p}\,\dd s>0\).
This proves \eqref{eq:seed}.
\end{proof}

\subsection{Subcritical case}

\begin{proof}[Proof of Theorem~\ref{thm:blowup} in the subcritical case \(\alpha_{\rm C}(p,q;3)>0\)]
Assume \(\alpha_{\rm C}(p,q;3)>0\), and suppose, for contradiction,
that \(T_m=\infty\).
For \(T\ge256\), consider the region
\begin{equation*}
        T/4\le t\le T/2,\qquad |x|\le\sqrt t/4.
\end{equation*}
Here
\begin{equation*}
        t\ge\frac T4\ge64,\qquad
        \frac tT\le\frac12,\qquad
        \frac{|x|}{\sqrt T}\le\frac1{4\sqrt2}<\frac12.
\end{equation*}
Hence \(\Psi_T=1\), and Lemma~\ref{lem:cone-seed} applies.
Consequently,
\begin{align*}
        J_T
        &\gtrsim
        \int_{T/4}^{T/2}t^{-q\min\{p/2,3/2\}}
        |B_{\sqrt t/4}|\,\dd t\\
        &\gtrsim T^{5/2-q\min\{p/2,3/2\}}.
\end{align*}
Combining this with \eqref{eq:test-upper} and using
\begin{equation*}
        1+\frac{q+1}{pq-1}
        -q\min\left\{\frac p2,\frac32\right\}
        =q\alpha_{\rm C}(p,q;3)
\end{equation*}
gives
\[
        T^{q\alpha_{\rm C}(p,q;3)}\lesssim 1.
\]
For the fixed profiles \(\phi,\psi\), this is impossible as \(T\to\infty\).
\end{proof}

\subsection{Critical case\texorpdfstring{: $1<p<3$}{: 1<p<3}}

We now assume
\begin{equation}\label{eq:critical-cone-relation-main}
        1<p<3,\qquad q=\frac{3p+2}{p^2}.
\end{equation}
In particular, \(pq=3+2/p\).  
We transfer the parabolic lower bound once to the wave component.
Note that Kirchhoff's formula gives
\begin{equation}\label{eq:wave-duhamel}
        L_W[F](t,x)
        =\frac1{4\pi}\int_0^t(t-s)
        \int_{\mathbb S^2}F(s,x+(t-s)\omega)\,\dd\omega\,\dd s.
\end{equation}

\begin{lemma}\label{lem:core-to-cone-critical}
The following estimate holds.
\begin{equation}\label{eq:core-to-cone-critical-main}
        v(t,x)\gtrsim
        t^{-1}(t-|x|)^{-1/p}
\end{equation}
whenever \(0<t<T_m\) and \(64\le t-|x|\le t/2\).
\end{lemma}

\begin{proof}
Put \(\zeta=t-|x|\) and restrict \eqref{eq:wave-duhamel} to
\[
        \zeta\le s\le\zeta+\frac{\sqrt\zeta}{16}.
\]
Since \(\zeta\ge64\), we have \(\sqrt\zeta/16\le\zeta/2\).
Writing \(\rho=t-s\), we obtain
\begin{align*}
        \zeta&\le s\le\frac{3\zeta}{2}\le\frac{3t}{4},\\
        \bigl||x|-\rho\bigr|&=|s-\zeta|
        \le\frac{\sqrt\zeta}{16}\le\frac{\sqrt s}{16},\\
        \frac t2&\le|x|\le t,\qquad
        \frac t4\le\rho\le t.
\end{align*}
The spherical reduction formula \eqref{eq:spherical-reduction-global} gives
\begin{align*}
        \left|\left\{\omega\in\mathbb S^2:
        |x+\rho\omega|\le\frac{\sqrt s}{4}\right\}\right|
        &=\frac{2\pi}{|x|\rho}\int_{||x|-\rho|}^{|x|+\rho}
   \mathbf 1_{[0,\sqrt{s}/4]}(\lambda)\,\lambda\, \dd\lambda\\
        &=\frac{2\pi}{|x|\rho}\int_{||x|-\rho|}^{\sqrt{s}/4}\lambda\,\dd\lambda\\
        &=\frac{\pi\bigl(s/16-(|x|-\rho)^2\bigr)}{|x|\rho}\\
        &\ge\frac{\pi(s/16-s/256)}{t^2}
        =\frac{15\pi s}{256t^2}.
\end{align*}
Here \(t\ge2\zeta\ge128\), and
\begin{equation*}
        \bigl||x|-\rho\bigr|\le\frac{\sqrt s}{16}
        <\frac{\sqrt s}{4}<\frac{\sqrt t}{4}
        <\frac t2\le|x|+\rho.
\end{equation*}
For every \(\omega\in\mathbb S^2\) satisfying
\(|x+\rho\omega|\le\sqrt{s}/4\), Lemma~\ref{lem:cone-seed} gives
\[
    u(s,x+\rho\omega)^q\gtrsim s^{-pq/2}.
\]
Since \(1-pq/2=-1/2-1/p\), it follows that
\begin{align*}
        v(t,x)&\gtrsim t^{-1}
        \int_\zeta^{\zeta+\sqrt\zeta/16}s^{-1/2-1/p}\,\dd s\\
        &\ge t^{-1}\frac{\sqrt\zeta}{16}
        \left(\frac{3\zeta}{2}\right)^{-1/2-1/p}\\
        &=\frac{1}{16}\left(\frac32\right)^{-1/2-1/p}
        t^{-1}\zeta^{-1/p}.
\end{align*}
This proves \eqref{eq:core-to-cone-critical-main}.
\end{proof}

\begin{proof}[Proof of Theorem~\ref{thm:blowup} in the critical case \(\alpha_{\rm C}(p,q;3)=0\) when \(1<p<3\)]
Suppose, for contradiction, that \(T_m=\infty\).
On \eqref{eq:critical-cone-relation-main}, \eqref{eq:test-upper} reads
\begin{equation}
    \label{eq:esti_upperIT}
    I_T\lesssim T^{(3-p)/2}.
\end{equation}
For \(T\ge2^{24}\), restrict the integral defining \(I_T\) to
\begin{equation*}
        \frac{\sqrt T}{8}\le t\le\frac{\sqrt T}{4},
        \qquad
        64\le\zeta:=t-|x|\le\frac t4.
\end{equation*}
Since \(t\ge\sqrt T/8\ge512\), the interval \([64,t/4]\) is nonempty.
Moreover,
\begin{equation*}
        \frac tT\le\frac1{4\sqrt T}\le\frac12,\qquad
        \frac{|x|}{\sqrt T}\le\frac t{\sqrt T}\le\frac14<\frac12.
\end{equation*}
Thus \(\Psi_T=1\).  Since \(64\le\zeta\le t/4\le t/2\) and
\(3t/4\le|x|=t-\zeta\le t\), Lemma~\ref{lem:core-to-cone-critical}
applies.  Also, \(\log T\ge24\log2\) gives
\begin{equation*}
        \log\frac{t}{256}
        \ge\frac12\log T-\log2048
        \ge\frac1{24}\log T.
\end{equation*}
Using polar coordinates, we obtain
\begin{align}
        I_T
        &\gtrsim
        \int_{\sqrt T/8}^{\sqrt T/4}
        t^{2-p}\int_{64}^{t/4}\frac{\dd\zeta}{\zeta}\,\dd t\nonumber\\
        &=
        \int_{\sqrt T/8}^{\sqrt T/4}
        t^{2-p}\log\frac{t}{256}\,\dd t\nonumber\\
        &\ge\frac{1}{24}\log T
        \int_{\sqrt T/8}^{\sqrt T/4}t^{2-p}\,\dd t\nonumber\\
        &=\frac{1}{24}\frac{4^{p-3}-8^{p-3}}{3-p}
        T^{(3-p)/2}\log T.
        \label{eq:esti_lowerIT}
\end{align}
Thus, \eqref{eq:esti_upperIT} and \eqref{eq:esti_lowerIT} give \(\log T\lesssim 1\), which is impossible as \(T\to\infty\).
\end{proof}

\subsection{Critical case\texorpdfstring{: $3\le p<5$}{: 3<=p<5}}

We finally consider
\begin{equation}\label{eq:heat-critical-relation}
        3\le p<5,\qquad q=\frac{2p+5}{3p}.
\end{equation}
Here the integral upper bounds themselves suffice.

\begin{proof}[Critical case \(3\le p<5\) of Theorem~\ref{thm:blowup}]
Suppose, for contradiction, that \(T_m=\infty\).
Since \((p+1)/(pq-1)=3/2\), \eqref{eq:test-upper} gives \(I_T\lesssim 1\).
For every fixed \((t,x)\), we have \(\Psi_T(t,x)=1\) whenever
\(T\ge\max\{1,2t,4|x|^2\}\), and hence
\(\Psi_T(t,x)\to1\) as \(T\to\infty\).
Fatou's lemma therefore yields
\begin{equation}\label{eq:heat-critical-integrability}
        \int_0^\infty\int_{\R^3}v(t,x)^p\,\dd x\,\dd t<\infty.
\end{equation}

The derivatives of \(\Psi_T\) are supported in
\[
        E_T:=\{(t,x)\in[0,T]\times B_{\sqrt T}:
        t\ge T/2\ \hbox{or}\ |x|\ge\sqrt T/2\}.
\]
Retaining this support in the proof of \eqref{eq:test-wave}, we obtain
\begin{equation}\label{eq:test-wave-tail}
        J_T\lesssim T^{-1+5/(2p')}
        \left(\int_{E_T}v(t,x)^p\,\dd x\,\dd t\right)^{1/p}.
\end{equation}
On \eqref{eq:heat-critical-relation},
\[
        -1+\frac5{2q'}
        +\frac1q\left(-1+\frac5{2p'}\right)=0.
\]
Moreover, \(\Psi_T(0,x)=1\) on \(\operatorname{supp}\phi\) for \(T\ge16\).
Combining \eqref{eq:test-damped} and \eqref{eq:test-wave-tail} gives
\[
        I_T+\int_{\R^3}\phi(x)\,\dd x
        \lesssim \left(\int_{E_T}v(t,x)^p\,\dd x\,\dd t\right)^{1/(pq)}
        \longrightarrow0.
\]
Indeed,
\begin{align*}
        0\le\int_{E_T}v(t,x)^p\,\dd x\,\dd t
        &\le\int_{T/2}^{\infty}\int_{\R^3}v(t,x)^p\,\dd x\,\dd t\\
        &\quad+\int_0^\infty\int_{|x|\ge\sqrt T/2}
        v(t,x)^p\,\dd x\,\dd t
        \longrightarrow0
\end{align*}
by \eqref{eq:heat-critical-integrability}.
This contradicts \(\int_{\R^3}\phi>0\).
\end{proof}

Together with the strict case, these two critical cases cover
\eqref{eq:blowup-index-condition} by Remark~\ref{rem:alphaCP-blowup-equivalence}.
This completes the proof of Theorem~\ref{thm:blowup}.

\appendix
\section{Proof of local well-posedness}\label{app:local-sobolev}

In this appendix, we prove Proposition~\ref{prop:local-sobolev}.  
We consider a slightly more general system
\begin{equation}\label{eq:app-Nakao-fg}
\left\{
\begin{aligned}
        &\partial_t^2 u-\Delta u+\partial_t u=f(v),&&\qquad t>0,~x\in\R^3,\\
        &\partial_t^2 v-\Delta v=g(u),&&\qquad t>0,~x\in\R^3,\\
        &u(0)=u_0,\quad \partial_t u(0)=u_1,\\
        &v(0)=v_0,\quad \partial_t v(0)=v_1,
\end{aligned}
\right.
\end{equation}
where \(f,g\in C^1(\R)\) and \(f(0)=g(0)=0\).  The case of Proposition~\ref{prop:local-sobolev} corresponds to \(f(s)=|s|^p\) and \(g(s)=|s|^q\).
The Duhamel formulation is
\begin{equation}\label{eq:app-sobolev-duhamel}
\left\{
\begin{aligned}
        u(t)&=\partial_tD(t)u_0+D(t)(u_0+u_1)+\int_0^tD(t-\tau)f(v(\tau))\,\dd\tau,\\
        v(t)&=\partial_t W(t)v_0+W(t)v_1+\int_0^tW(t-\tau)g(u(\tau))\,\dd\tau.
\end{aligned}
\right.
\end{equation}
Here \(D(t)\) is the propagator of the damped wave equation defined in \eqref{eq:D-decomposition}, and \(W(t)\) is the propagator of the wave equation in \eqref{eq:W-explicit}.

We use the following standard linear estimates on finite time intervals.

\begin{lemma}\label{lem:app-sobolev-linear}
There exists a constant \(C_0>0\) such that, for every \(T>0\), the following holds for \(A=D,W\) and \(0\le t\le T\) with \(C_T:=C_0(1+T)\):
\begin{align*}
        \norm{K_A[\phi,\psi](t)}_{H^2}
        +\norm{\partial_tK_A[\phi,\psi](t)}_{H^1}
        \le
        C_T\left(\norm{\phi}_{H^2}+\norm{\psi}_{H^1}\right),
\end{align*}
for all \((\phi,\psi)\in H^2(\R^3)\times H^1(\R^3)\).
Moreover, for every \(F\in L^1(0,T;L^2(\R^3))\),
\begin{align*}
    \left\|\int_0^tA(t-s)F(s)\,\dd s\right\|_{L^2}
        \le C_T\int_0^t\norm{F(s)}_{L^2}\,\dd s.
\end{align*}
In addition, for every \(F\in L^1(0,T;H^1(\R^3))\),
\begin{align*}
    \left\|\int_0^t A(t-s)F(s)\,\dd s\right\|_{H^2}
        +
        \left\|\partial_t\int_0^t A(t-s)F(s)\,\dd s\right\|_{H^1}
        \le
        C_T\int_0^t \norm{F(s)}_{H^1}\,\dd s.
\end{align*}
For such \(F\) and \((\phi,\psi)\in H^2(\R^3)\times H^1(\R^3)\), \(K_A[\phi,\psi]+L_A[F]\) belong to \(C([0,T];H^2(\R^3))\cap C^1([0,T];H^1(\R^3))\). 
If, in addition, \(F\in C([0,T];L^2(\R^3))\), then it belongs to \(C^2([0,T];L^2(\R^3))\).
\end{lemma}

The estimates of the damped wave equation follow from the \(L^2\) estimates in Matsumura~\cite[Section~2.1, Lemma~1]{Matsumura} with \(m=2\), by a density argument and Minkowski's inequality; time continuity follows from the multiplier formulas.
For the wave equation, see \cite[Proposition~3.19, Theorem~3.24, and Proposition~3.30]{Katayama2017}; the extension to \(L^1\)-in-time forcing follows by a density argument.
For either equation, the \(C^2([0,T];L^2)\) regularity follows from the equation itself.

For simplicity, we here use the following norms on \(H^m(\R^3)\), \(m=1,2\):
\[
    \norm{f}_{H^m}:=\norm{(1+|\xi|^2)^{m/2}\widehat f(\xi)}_{L^2_\xi},
    \qquad
    \widehat f(\xi):=(2\pi)^{-3/2}\int_{\R^3}e^{-ix\cdot\xi}f(x)\dd x,
\]
where the Fourier transform is initially defined on Schwartz functions and extended to \(L^2\).
Then, the inequality \(\|f\|_{L^\infty}\le \|f\|_{H^2}\) holds for \(f\in H^2\).
Indeed, the Fourier inversion theorem and the Cauchy--Schwarz inequality give
\[
    |f(x)|
    \le (2\pi)^{-3/2}
    \left(\int_{\R^3}\frac{\dd\xi}{(1+|\xi|^2)^2}\right)^{1/2}\norm{f}_{H^2}
    =\frac{1}{\sqrt{8\pi}}\norm{f}_{H^2},
\]
since \(4\pi\int_0^\infty r^2(1+r^2)^{-2}\dd r=\pi^2\).

\begin{proof}[Proof of Proposition~\ref{prop:local-sobolev}]
We prove the result for \eqref{eq:app-Nakao-fg}.  For \(0<T<1\), Lemma~\ref{lem:app-sobolev-linear} gives
\begin{align*}
    &\norm{K_D[u_0,u_1]}_{L^\infty(0,T;H^2)}
     +\norm{K_W[v_0,v_1]}_{L^\infty(0,T;H^2)}\\
    &\qquad\le 2C_0\bigl(\norm{u_0}_{H^2}+\norm{u_1}_{H^1}
                        +\norm{v_0}_{H^2}+\norm{v_1}_{H^1}\bigr).
\end{align*}
Choose
\[
    \rho:=\max\bigl\{1,\,4C_0\bigl(\norm{u_0}_{H^2}+\norm{u_1}_{H^1}
                        +\norm{v_0}_{H^2}+\norm{v_1}_{H^1}\bigr)\bigr\}
\]
and set
\begin{equation*}
        \mathcal B_T(\rho):=\left\{(u,v)\in L^\infty(0,T;H^2)^2:
        \norm{u}_{L^\infty(0,T;H^2)}+\norm{v}_{L^\infty(0,T;H^2)}\le\rho\right\}.
\end{equation*}
We equip this set with the metric
\begin{equation*}
        d((u^1,v^1),(u^2,v^2))
        :=\norm{u^1-u^2}_{L^\infty(0,T;L^2)}
        +\norm{v^1-v^2}_{L^\infty(0,T;L^2)}.
\end{equation*}
The metric space \((\mathcal B_T(\rho),d)\) is complete.  
Indeed, every Cauchy sequence with respect to this distance \(d\) has a limit in \(L^\infty(0,T;L^2)^2\).  
Since the sequence is bounded in \(L^\infty(0,T;H^2)^2\), weak-star compactness and lower semicontinuity show that the limit belongs to \(\mathcal B_T(\rho)\).

For \(R>0\), define
\begin{equation*}
        M(R):=\max\left\{\sup_{|s|\le R}|f'(s)|,
        \sup_{|s|\le R}|g'(s)|\right\}.
\end{equation*}
If \((u,v)\in\mathcal B_T(\rho)\), then \(\norm{u(t)}_{L^\infty}+\norm{v(t)}_{L^\infty}\le\rho\) for a.e. \(t\).  Since \(f(0)=g(0)=0\), the mean value formula and the chain rule give
\begin{equation}\label{eq:app-fg-H1-bound}
        \norm{f(v)}_{L^\infty(0,T;H^1)}
        +\norm{g(u)}_{L^\infty(0,T;H^1)}
        \le M(\rho)\left(\norm{v}_{L^\infty(0,T;H^1)}
                        +\norm{u}_{L^\infty(0,T;H^1)}\right).
\end{equation}
Moreover, for \((u^i,v^i)\in\mathcal B_T(\rho)\), \(i=1,2\),
\begin{equation}\label{eq:app-fg-L2-Lip}
\begin{aligned}
&\norm{f(v^1(t))-f(v^2(t))}_{L^2}
  +\norm{g(u^1(t))-g(u^2(t))}_{L^2}\\
&\qquad\le M(\rho)
\left(\norm{v^1(t)-v^2(t)}_{L^2}
+\norm{u^1(t)-u^2(t)}_{L^2}\right)
\end{aligned}
\end{equation}
for a.e. \(t\in(0,T)\).

Define \(\Phi:\mathcal B_T(\rho)\to L^\infty(0,T;H^2)^2\) by the right-hand side of \eqref{eq:app-sobolev-duhamel}.
Here and below, the norm of a pair is the sum of the component norms.
For \((u,v)\in\mathcal B_T(\rho)\), using \eqref{eq:app-fg-H1-bound} and Lemma~\ref{lem:app-sobolev-linear}, we obtain
\begin{align*}
\norm{\Phi(u,v)}_{L^\infty(0,T;H^2)^2}
&\le \frac{\rho}{2}+2C_0M(\rho)T
\left(\norm{u}_{L^\infty(0,T;H^1)}
+\norm{v}_{L^\infty(0,T;H^1)}\right)\\[2mm]
&\le \frac{\rho}{2}+2C_0M(\rho)T\rho .
\end{align*}
Choose \(0<T<1\) so that \(2C_0M(\rho)T\le 1/2\).
Then \(\Phi\) maps \(\mathcal B_T(\rho)\) into itself.
The \(L^2\)-estimates in Lemma~\ref{lem:app-sobolev-linear}, together with \eqref{eq:app-fg-L2-Lip}, give
\begin{equation*}
        d(\Phi(u^1,v^1),\Phi(u^2,v^2))
        \le 2C_0M(\rho)T\,d((u^1,v^1),(u^2,v^2))
\end{equation*}
for every \((u^1,v^1), (u^2,v^2)\in \mathcal B_T(\rho)\).
Thus \(\Phi:\mathcal B_T(\rho) \to \mathcal B_T(\rho)\) is a contraction.
Hence there is a unique \((u,v)\in\mathcal B_T(\rho)\) such that \((u,v)=\Phi(u,v)\), that is, \((u,v)\in L^\infty(0,T;H^2)^2\) satisfies \eqref{eq:app-sobolev-duhamel}.

By \eqref{eq:app-fg-H1-bound}, the nonlinear terms belong to \(L^\infty(0,T;H^1)\subset L^1(0,T;H^1)\).  Lemma~\ref{lem:app-sobolev-linear} gives
\begin{equation*}
        u,v\in C([0,T];H^2)\cap C^1([0,T];H^1).
\end{equation*}
Since \(u,v\in C([0,T];H^2)\), the maps \(t\mapsto f(v(t))\) and \(t\mapsto g(u(t))\) belong to \(C([0,T];L^2)\).  The equations then imply \(\partial_t^2 u,\partial_t^2 v\in C([0,T];L^2)\).  This proves the regularity \eqref{eq:local-sobolev-regularity-main} on \([0,T]\).

The same construction starts from any time \(t_0\) with data \((u(t_0),\partial_t u(t_0))\) and \((v(t_0),\partial_t v(t_0))\) in \(H^2\times H^1\).  Hence one obtains a maximal solution on \([0,T_m)\).  Uniqueness on overlapping intervals follows from the same contraction argument on sufficiently short subintervals.

We first prove the standard blow-up alternative of Sobolev norms:
\begin{equation}\label{eq:local-sobolev-alt-main}
        T_m<\infty
        \quad\Longrightarrow\quad
        \lim_{t\uparrow T_m}
        \left(
        \norm{(u,v)(t)}_{H^2}
        +\norm{(\partial_t u,\partial_t v)(t)}_{H^1}
        \right)=\infty .
\end{equation}
Indeed, if \eqref{eq:local-sobolev-alt-main} failed, then there would be \(K>0\) and \(t_n\uparrow T_m\) such that the \(H^2\times H^1\)-norm of the solution at \(t_n\) is bounded by \(K\).  The local existence time obtained above depends only on such a bound.  Therefore the solution can be restarted at \(t_n\) for a time \(\delta(K)>0\) independent of \(n\).  Taking \(n\) large enough gives an extension beyond \(T_m\), contradicting maximality.  This proves \eqref{eq:local-sobolev-alt-main}.

We next turn the alternative \eqref{eq:local-sobolev-alt-main} into the \(L^\infty\)-alternative stated in Proposition~\ref{prop:local-sobolev}.  
Set
\begin{equation*}
        E(t):=\norm{u(t)}_{H^2}+\norm{\partial_t u(t)}_{H^1}
        +\norm{v(t)}_{H^2}+\norm{\partial_t v(t)}_{H^1}.
\end{equation*}
Suppose that \(T_m<\infty\) and that the limsup in \eqref{eq:local-sobolev-Linfty-alt-main} is finite.  
Then, there exists \(K>0\) such that
\begin{equation*}
        \sup_{0\le t<T_m}
        \left(\norm{u(t)}_{L^\infty}+\norm{v(t)}_{L^\infty}\right)\le K .
\end{equation*}
Since \(f(0)=g(0)=0\), the mean value formula and the chain rule give
\begin{equation*}
        \norm{f(v(t))}_{H^1}+\norm{g(u(t))}_{H^1}
        \le M(K)\left(\norm{v(t)}_{H^1}+\norm{u(t)}_{H^1}\right)
        \le M(K)E(t)
        \qquad (0\le t<T_m).
\end{equation*}
Applying Lemma~\ref{lem:app-sobolev-linear} yields
\begin{equation*}
        E(t)\le C_{T_m}E(0)
        +C_{T_m}M(K)\int_0^t E(\tau)\,\dd\tau
        \qquad (0\le t<T_m).
\end{equation*}
Therefore, Gronwall's inequality gives
\begin{equation*}
        \sup_{0\le t<T_m}E(t)
        \le C_{T_m}E(0)e^{C_{T_m}M(K)T_m}
        <\infty .
\end{equation*}
This contradicts the alternative \eqref{eq:local-sobolev-alt-main}.  
Hence \eqref{eq:local-sobolev-Linfty-alt-main} holds.
\end{proof}

\begin{proof}[Proof of Corollary~\ref{cor:compatibility}]
By uniqueness of bounded mild solutions on finite time intervals, the maximal solution of Proposition~\ref{prop:local-sobolev} satisfies \eqref{eq:global-bound} on \([0,T_m)\).
Hence the \(L^\infty\) blow-up alternative \eqref{eq:local-sobolev-Linfty-alt-main} implies \(T_m=\infty\).
\end{proof}

\medskip
\noindent\textbf{Statements and Declarations.}\quad
The authors declare that they have no competing interests.

\medskip
\noindent\textbf{Acknowledgments.}\quad
V. G. is partially supported by the Gruppo Nazionale per l'Analisi Matematica and by the INdAM-GNAMPA Project, cod. CUP E53C25002010001.
K. K. is partially supported by JSPS KAKENHI Grant Numbers JP24K16955 and JP24KJ0010.

\end{document}